\documentclass[a4paper,12pt]{article}
\usepackage[english]{babel}
\usepackage{amsfonts,amssymb,amsmath}
\begin{document}

 \centerline{\bf\Large Large-gradient and Lagrange singularities}
 \vskip 2mm  
 \centerline{\bf\Large of solutions of a quasilinear parabolic equation}

 \vskip 7mm 
 \centerline{\bf S.V.~Zakharov} 

 \vskip 5mm 

\begin{center}
Institute of Mathematics and Mechanics,\\
Ural Branch of the Russian Academy of Sciences,\\
16, S.Kovalevskaja street, 620990, Ekaterinburg, Russia
\end{center}

\vskip 7mm 

\textbf{Abstract.}
Asymptotic solutions of 
a quasilinear parabolic equation
with a small parameter at the higher derivative
are constructed near large-gradient and Lagrange singularities
of $A$-type, which represent interest for studying
processes of shock waves formation
in physical media with a small nonzero viscosity. 
 
Key words: parabolic equation, Burgers equation, Cole--Hopf transform, asymptotics, large gradient, Lagrange singularities, Whitney fold. 

Mathematics Subject Classification: 35K15, 35K59.
 
 \vskip 15mm 
 \section{Introduction} 

A simplest model of the motion of continuum,
which takes into account nonlinear effects and dissipation,
 is the equation of nonlinear diffusion $$ 
 \frac{\partial u}{\partial t} + 
 u\frac{\partial u}{\partial x} = 
 \varepsilon \frac{\partial^2 u}{\partial
x^2} 
 $$ 
for the first time presented by
H.~Bateman~\cite{Bate} and J.~Burgers~\cite{bu}.
 This equation  is used in studying
 the evolution of a wide class of physical systems and probabilistic
process, acoustic waves in fluid and gas~\cite{wh}.

 In the present paper, results of investigations
of the asymptotic  behavior of solutions near singular points
in the Cauchy problem for the more general
quasilinear parabolic equation
 \begin{eqnarray} 
 \label{main} 
 \frac{\partial u}{\partial t} + 
 \frac{\partial \varphi(u)}{\partial
x} = 
 \varepsilon \frac{\partial^2 u}{\partial
x^2}, 
 & \quad t\geqslant t_0, \\ 
 \label{ic}\phantom{\frac{1}{1}} 
 u(x,t_0) = q(x), & \quad x\in\mathbb{R},
 \end{eqnarray} 
are given.
 We assume that $\varepsilon>0,$ 
the function $\varphi$ is infinitely differentiable
 and its second derivative is strictly positive. 
 The initial function~$q$ is bounded and smooth.

The interest to study the behavior of solutions near
singular points is explained, in
particular, the fact that such singular events
 last very  little time and, however, they determine
 all subsequent behavior of the system in many respects. 
The asymptotic behavior of solutions in neighborhoods
of singular points is directly connected with
constructing approximations in neighborhoods
of shock waves,
which is important for applied problems
in the mechanics of continua. 
 
 Although the types of singular points of
 solutions are in detail classified~\cite{as} 
 and processes of shock waves formation
 for the degenerate first-order equation are studied~\cite{bb}, 
 constructing asymptotic series in the viscosity 
parameter $\varepsilon$ for an equation of a more general form~(\ref{main}) 
 in every specific case is a separate and very complex problem. 
 For example, in the problem about
 the transition of a weak discontinuity
 into a strong one~\cite{gc} 
there is still an open question of matching
asymptotic series. 

\vskip 10mm
\section{Singularity generated by a large initial gradient} 
 
First, we consider a singular point of the solution
with two small parameters~\cite{2ps}
 in the case, when the initial condition has the
form $$ 
 \phantom{\frac{1}{1}} 
 u(x,0,\varepsilon,\rho) = \nu ( {x}{\rho}^{-1}), 
 \quad x\in\mathbb{R}, \quad \rho>0, 
 $$ 
 where function $\nu$ is infinitely differentiable and bounded,
 and $\rho$ is the second small parameter. 
 It is proved that in this case
under conditions $$ 
 \nu(\sigma) = \sum\limits_{n=0}^{\infty} 
 \frac{\nu^{\pm}_n}{\sigma^{n}}, 
 \qquad 
 \sigma\to \pm\infty, 
 \qquad 
 (\nu^-_0>\nu^+_0) 
 $$ 
 for the solution of problem $(\ref{main})$--$(\ref{ic})$ 
 as
$$
\varepsilon\to 0
\qquad
\mbox{and}
\qquad \mu=\rho/\varepsilon\to 0
$$ 
 in the strip $$ 
 \{ (x,t) : x\in\mathbb{R},\ 0 \leqslant
t\leqslant T\} 
 $$ 
 there holds the asymptotic formula $$ 
 {u}(x,t,\varepsilon,\rho)= 
 h_0\left(  \frac{x}{\rho}, 
 \frac{\varepsilon t}{\rho^2}\right)
- R_{0,0,0}\left( 
 \frac{x}{2\sqrt{\varepsilon t}} 
 \right)+ 
 \Gamma\left( 
 \frac{x}{\varepsilon}, 
 \frac{t}{\varepsilon} 
 \right) 
 +O\left(\mu^{1/2}\ln\mu\right), 
 $$ 
 where $$ 
 h_0(\sigma,\omega) = 
 \frac{1}{2\sqrt{\pi\omega}} 
 \int\limits_{-\infty}^{\infty} 
 \nu(s) \exp\left[ -\frac{(\sigma-s)^2}{4\omega}\right] 
 ds, 
 $$ 
 $$ 
 R_{0,0,0}(z) = \nu^{-}_{0} 
 \mathrm{erfc}(z) + \nu^{+}_{0} \mathrm{erfc}(-z), 
 $$ 
 $$ 
 \mathrm{erfc} (z) = \frac{1}{\sqrt{\pi}} 
 \int\limits_{z}^{+\infty} \exp(-y^2)\,
dy, 
$$
$$
\sigma = \frac{x}{\rho},
 \qquad 
\omega = \frac{\varepsilon t}{\rho^2},
 \qquad 
 z = \frac{\sigma}{2\sqrt{\omega}}, 
 $$ 
the function $\Gamma$ is the solution 
of the equation in the inner variables 
$$
\eta = x/\varepsilon,
\qquad \theta = t/\varepsilon,
$$
 $$ 
 \frac{\partial \Gamma}{\partial \theta} + 
 \frac{\partial \varphi(\Gamma)}{\partial \eta} - \frac{\partial^2 \Gamma}{\partial \eta^2} = 0 
 $$ 
 with the initial condition $$ 
 \Gamma(\eta,0) = 
 \begin{cases} 
 \nu^-_0, & \eta<0, \\ 
 \nu^+_0, & \eta > 0. 
 \end{cases} 
 $$ 

 In addition in~\cite{zz},
using the renormalization method~\cite{gv}
the following asymptotic  formula 
is obtained: $$ 
 {u}(x,t,\varepsilon,\rho)= 
 \frac{1}{\nu^+_0 - \nu^-_0} 
 \int\limits_{-\infty}^{\infty} 
 \Gamma\left( \frac{x-\rho s}{\varepsilon}, \frac{t}{\varepsilon}\right) 
 \nu'(s)\, ds 
 +O\left(\mu^{1/4}\right). 
 $$

\vskip 10mm 
\section{Fold singularity} 

 Consider the Cauchy problem~(\ref{main})--(\ref{ic})
 for the quasilinear parabolic equation 
 when the solution of the limit problem $(\varepsilon=0)$ 
 has a point of gradient catastrophe. 
 A.M. Il'in studied this problem  
in the case when in the strip $$\{(x,t):~t_0\leqslant t\leqslant T,~x\in\mathbb{R}\}$$ 
the limit solution is a smooth  function everywhere except for one
smooth discontinuity curve $$ 
 \{ (x,t) : x= s(t),\ t\geqslant t^{*}>
t_0\}. 
 $$ 
A detailed presentation of his results can be found in
monograph~\cite{ib}, 
 where the asymptotics of the solution
as $\varepsilon\to 0$ is constructed and justified with
an arbitrary accuracy. 
Under a suitable choice of independent variables,
the  singular point $(s(t^{*}),t^{*})$ 
coincides with the origin and
 in its neighborhood 
the following stretched variables are introduced:
 $$ 
 \xi = \varepsilon^{-3/4} x, 
 \quad 
 \tau = \varepsilon^{-1/2} t.
 $$ 
An expansion of the solution is sought in the form
of the series
 $$ 
 w = \sum\limits_{k=1}^{\infty} \varepsilon^{k/4} 
 \sum\limits_{j=0}^{k-1} w_{k,j}(\xi,\tau) 
 \ln^{j} \varepsilon^{1/4}.
 $$ 
Substituting it into equation~(\ref{main}), 
for coefficients $w_{k,j}$ we obtain
the recurrence system $$ 
 \frac{\partial w_{1,0}}{\partial \tau} 
 + \varphi''(0) w_{1,0}\frac{\partial w_{1,0}}{\partial \xi} 
 -\frac{\partial^2 w_{1,0}}{\partial \xi^2} =0, 
 $$ 
 $$ 
 \frac{\partial w_{k,j}}{\partial \tau} 
 + \varphi''(0) \frac{(\partial w_{1,0} w_{k,j})}{\partial \xi} 
 -\frac{\partial^2 w_{k,j}}{\partial \xi^2} =E_{k,j}, 
 \quad k\geqslant 2.
 $$ 
 These equations should be supplied with the conditions $$ 
 w_{k,j}(\xi,\tau) =W_{k,j}(\xi,\tau), 
 \quad \tau\to -\infty, 
 $$ 
 where $W_{k,j}(\xi,\tau)$ is the sum of all
coefficients at $\varepsilon^{k/4}\ln^{j} \varepsilon^{1/4}$
 in the reexpansion of the asymptotics far from
the singularity (the outer expansion) 
 in terms of the inner variables. 
 
 Investigation of solutions of this system is
the central and most laborious task in this problem. 
 It is proved that there exist solutions $w_{k,j}(\xi,\tau)$ 
 for $k\geqslant 2$, $0 \leqslant j \leqslant
k-1$, 
which are infinitely differentiable for all $\xi$ and $\tau$. 
 
 Observe separately the properties of the leading term,
 which is found with the help of the Cole--Hopf transform $$ 
 w_{1,0} (\xi,\tau) = - \frac{2}{\varphi''(0)\Lambda(\xi,\tau)}
 \frac{\partial \Lambda(\xi,\tau)}{\partial \xi}, 
 $$ 
 where $$ 
 \Lambda(\xi,\tau) =\int\limits_{-\infty}^{\infty} 
 \exp\left( -2 s^4 + \tau s^2 - \xi s  \right) ds
 $$ 
 is a real-valued analog of the Pearcey function~\cite{tp}, 
 which is a solution of the heat equation. 
The argument of the exponent is the generating family
of the Lagrange singularity $A_3$, 
 see~\cite{as}. 
 
 \vskip 5mm 
 \noindent 
 \textbf{Theorem~1.} 
 {\it The function $w_{1,0}$ 
 satisfies the asymptotic relations $$ 
 w_{1,0} (\xi,\tau) = 
 [\varphi''(0)]^{-1} H(\xi,\tau) + 
 \sum\limits_{l=1}^{\infty} h_{1-4l}(\xi,\tau), 
 \qquad 
 3[H(\xi,\tau)]^2-\tau\to \infty, 
 $$ 
 $$ 
 (\xi,\tau)\in \Omega_1=\mathbb{R}^2 \setminus 
 \{ |\xi|< \tau^{\gamma_1-1/2}, \ \tau >0, \ 0< \gamma_1< 2\}, 
 $$ 
 where $H(\xi,\tau)$ is the Whitney fold
function, $$ H^3 - \tau H + \xi = 0, $$ 
 $h_l(\xi,\tau)$ are homogeneous functions
of power~$l$ 
 in $H(\xi,\tau)$, $\sqrt{-\tau}$
and $\sqrt{3[H(\xi,\tau)]^2-\tau}$, 
 which are polynomials in $H(\xi,\tau)$, $\tau$
and $(3[H(\xi,\tau)]^2-\tau)^{-1}$, 
 $$ 
 w_{1,0}(\xi,\tau) = \sqrt{\tau} 
 \left( -\frac{\mathrm{th}\, z}{\varphi''(0)} 
 +\sum\limits_{k=1}^{\infty} \tau^{-2k}
q_k(z) 
 \right), 
 \qquad 
 \tau \to+\infty, 
 $$ 
 $$ 
 (\xi,\tau)\in \Omega_2=\{ |\xi|\tau^{1/2}< \tau^{\gamma_2}, 
 \ \gamma_1< \gamma_2 < 2\}, 
 $$ 
 where $z= \xi\sqrt{\tau}/2$, 
 $|q_k(z)|\leqslant M_k(1+|z|^k)$, 
 $q_k\in C^{\infty}(\mathbb{R}^1)$ are solutions
 of a recurrence system of ordinary differential
equations. 
 } 
 
The proof of the theorem is based
 on calculating the asymptotics
of the integral $\Lambda(\xi,\tau)$
by Laplace's method. 
 In the domain $\Omega_1$ the contribution
 into the asymptotics 
is given by one local maximum
 and in the domain $\Omega_2$
by two local maxima.

\vskip 10mm 
\section{Lagrange singularity $A_{2n+1}$} 

 Il'in investigated the problem
 with such an initial function $q(x)$ that 
$$ 
 \varphi'(q(x))= -x + x^{3}+ O(x^{4}), 
 \qquad x\to 0. 
 $$ 
In paper \cite{faa}, the more general condition 
$$ 
 \varphi'(q(x)) = -x + x^{2n+1}+ O(x^{2n+2}), 
 \qquad x\to 0, 
 $$ 
with an arbitrary natural~$n$
has been considered. 
 Then in the Cole--Hopf transform
 instead of the Pearcey function one should consider
the integral $$ 
 \int\limits_{-\infty}^{\infty} 
 \exp\left( - a s^{2n+2} +\theta s^2 -\eta s \right) ds
 $$ 
 corresponding a section 
of the Lagrange singularity $A_{2n+1}$. 
 
 Let us clear out scales of the inner variables,
which are introduced using the change 
$$ 
 x = \eta \varepsilon^{\sigma}, 
 \qquad 
 t = \theta \varepsilon^{\mu}, 
 \qquad u \sim \varepsilon^{\varkappa}. 
 $$ 
 
Since all terms in equation~(\ref{main}) 
should be of the same order,
we obtain the relation 
$$
 -\mu = \varkappa -\sigma = 1- 2\sigma. 
$$
 From the characteristic equation 
$$ x = y + (t + 1)\varphi'(q(y)) $$ 
 we obtain another relation
 $$ \sigma = \varkappa + \mu = (2n+1)\varkappa. $$ 
 From these relations we find 
\begin{equation}\label{smk} 
 \sigma = \frac{2n+1}{2n+2}, 
 \qquad 
 \mu = \frac{n}{n+1}, 
 \qquad 
 \varkappa = \frac{1}{2n+2}. 
 \end{equation} 
 
 Since $u = O (\varepsilon^{1/(2n+2)})$, 
 equation~(\ref{main}) becomes the Burgers equation,
 whose solution can be written 
in the form of the Cole--Hopf transform.
 Moreover, the coefficient at $s^{2n+2}$ 
 is determined from the condition 
of  matching the inner asymptotics and the outer
expansion $$ 
 u_{\mathrm{out}} \sim 
 \frac{ \varepsilon^{\varkappa} U_0(\eta,\theta)}{\varphi''(0)}, 
 \qquad 
 U_0^{2n+1} - \theta U_0 + \eta = 0. 
 $$ 
 
 \vskip 5mm 
 \noindent 
 \textbf{Theorem~2.} 
 {\it 
 In the domain $\Omega_{\varepsilon} =\{ (x,t)\, : 
 \, |x\varepsilon^{-\varkappa}| + |t| < K\varepsilon^{\mu},\ K>0 \}$ 
 for any $n\geqslant 2$
the~function $$ 
 u_{\mathrm{in}}(x,t,\varepsilon) = - 2 \varepsilon 
 [\varphi''(0) V(x,t,\varepsilon)]^{-1} 
 \frac{\partial V(x,t,\varepsilon)}{\partial
x}, 
 $$ 
 where $$ 
 V(x,t,\varepsilon) = \int\limits_{-\infty}^{\infty} 
 \exp\left( - \frac{2^{2n} s^{2n+2}}{n+1} + 
 \frac{t s^2}{\varepsilon^{\mu}} 
 - \frac{x s}{\varepsilon^{\sigma}} \right)
ds, 
 $$ 
 $$ 
 \sigma = \frac{2n+1}{2n+2}, 
 \qquad 
 \mu = \frac{n}{n+1}, 
 $$ 
 is an asymptotic solution of equation $(\ref{main})$ in
the following sense: 
 \begin{equation}\label{os} 
 \frac{\displaystyle\frac{\partial u_{\mathrm{in}}}{\partial t} + 
 	\frac{\partial \varphi(u_{\mathrm{in}})}{\partial
x} - \varepsilon \frac{\partial^2
u_{\mathrm{in}}}{\partial x^2} } 
 {\sup\limits_{(x,t)\in\Omega_{\varepsilon}} 
 \left\{ 
 \displaystyle 
 \left| \frac{\partial u_{\mathrm{in}}}{\partial
t} \right| 
 +\left| \frac{\partial \varphi(u_{\mathrm{in}})}{\partial x} \right| 
 +\left| \varepsilon \frac{\partial^2
u_{\mathrm{in}} }{\partial x^2} \right| 
 \right\} } 
 = O\left( \varepsilon^{\varkappa} \right), 
 \end{equation} 
where }
$$
 \varkappa = \frac{1}{2n+2}. 
$$ 
  
 \vskip 2mm 
 \noindent 
 \textbf{Proof.} 
 From the equation
$$
\frac{\partial V}{\partial t} 
=  \varepsilon \frac{\partial^2 V}{\partial x^2}
$$ 
 it follows that 
$$ 
 	\frac{\partial u_{\mathrm{in}}}{\partial t} + 
 \varphi''(0) u_{\mathrm{in}}\frac{\partial u_{\mathrm{in}}}{\partial x} = 
 \varepsilon \frac{\partial^2 u_{\mathrm{in}}}{\partial x^2}. 
 $$ 
 Then, using the relation 
$$ 
 \varphi'(u) = \varphi''(0) u + O(u^2),  
$$ 
 we obtain $$ 
 	\frac{\partial u_{\mathrm{in}}}{\partial t} + 
 	\frac{\partial \varphi(u_{\mathrm{in}})}{\partial x}
 - \varepsilon \frac{\partial^2
u_{\mathrm{in}}}{\partial x^2} 
 = \varphi'(u_{\mathrm{in}})\frac{\partial
u_{\mathrm{in}}}{\partial x} 
 -\varphi''(0) u_{\mathrm{in}}\frac{\partial
u_{\mathrm{in}}}{\partial x} 
 = O\left(u_{\mathrm{in}}^2 \frac{\partial
u_{\mathrm{in}}}{\partial x} \right). 
 $$ 
 Since $$ 
 u^2_{\mathrm{in}} = 4 \varepsilon^{2-2\sigma} 
 [\varphi''(0) V]^{-2} (V_{\eta})^2, 
 $$ 
 $$ 
 \frac{\partial u_{\mathrm{in}}}{\partial
x} = -2 \varepsilon^{1-2\sigma} 
 [\varphi''(0)]^{-1} 
 \left(\frac{V_{\eta\eta}}{V} 
 -\frac{(V_{\eta})^2}{V^2} \right), 
 $$ 
 we find $$ 
\frac{\partial u_{\mathrm{in}}}{\partial t} + 
 	\frac{\partial \varphi(u_{\mathrm{in}})}{\partial x}
 - \varepsilon \frac{\partial^2
u_{\mathrm{in}}}{\partial x^2} 
 = O(\varepsilon^{3-4\sigma} ). 
 $$ 
 Taking into account~(\ref{smk}) 
 and the order $O(\varepsilon^{2-3\sigma} )$
 of derivatives on the left-hand side, we obtain
estimate~(\ref{os}).


\vskip 10mm 
\begin{thebibliography}{99} 

\bibitem{Bate}
H. Bateman, 
{Some recent researches on the motion of fluids},
\textit{Monthly Weather Review},
\textbf{43}, (1915) 163-170.

\bibitem{bu}
J. Burgers,
 \textit{A Mathematical Model Illustrating the Theory of Turbulence},
 Advances in Applied Mechanics, Academic
Press, New York, 1948.

\bibitem{wh}
G.B. Whitham,
 \textit{Linear and Non-Linear Waves}, 
Wiley-Interscience, New York, 1974.

\bibitem{as}
V.I. Arnold,
{\it  Singularities of Caustics and Wave Fronts};
Mathematics and Its Applications Series, Vol. 62, Kluwer,
Dordrecht, 1990.

\bibitem{bb}
I.A. Bogaevsky, Metamorphoses of singularities
of minimum functions, and bifurcations of shock
waves of the Burgers equation with vanishing viscosity.
{\it Algebra i Analiz}, \textbf{1(4)}, (1989) 1-16.
 English translation in St. Petersburg (Leningrad) Math.
J., \textbf{1}, 4, (1990) 807-823.

\bibitem{gc}
   {S.V. Zakharov,}
Asymptotic solution of a Cauchy problem in a neighbourhood of a gradient catastrophe 
{\it Sbornik: Mathematics},
 \textbf{197}, 6, (2006) 835-851.

\bibitem{2ps}
   {S. V. Zakharov,}
 Two-parameter asymptotics in the Cauchy problem for a quasi-linear parabolic equation,
{\it  Asymptotic Analysis},
  \textbf{63}, 1-2, (2009) 49-54.

\bibitem{zz}
   {S. V. Zakharov,}
Cauchy problem for a quasilinear parabolic equation with a large initial gradient and low viscosity,
{\it Comp. Math. Math. Phys.},
 \textbf{50} Issue 4, (2010), 665-672. 

\bibitem{gv}
 N. Goldenfeld, J. Veysey,
Simple viscous flows: From boundary layers to the renormalization group,
\textit{Rev. Mod. Phys.}, \textbf{79}, 3, (2007) 883-927. 

\bibitem{ib}
A.M. Il'in,
{\it Matching of Asymptotic Expansions of Solutions
of Boundary Value Problems},
Am. Math. Soc., Providence, R.I., 1991.


\bibitem{tp}
T. Pearcey,
The Structure of an Electromagnetic Field in a Neighborhood of a Cusp of a Caustic, 
\textit{Phil. Mag.} 
\textbf{37}, (1946) 311. 

\bibitem{faa}
{S. V. Zakharov,}
Singularities of $A$ and $B$ types in asymptotic analysis of solutions of a parabolic equation,
\textit{Functional Analysis and Its Applications},
\textbf{49}, 4, (2015) 307-310. 

\end{thebibliography}
\end{document}